\documentclass{article}

\usepackage{amsmath,amssymb,amsthm}
\usepackage[latin1]{inputenc}

\newcommand{\dd}{\mathrm{d}}
\newcommand{\T}{^\mathrm{T}}
\newcommand{\xupt}{t,\mathbf{x}(t),\mathbf{u}(t),\psi_0,\boldsymbol{\psi}(t)}
\newcommand{\xup}{t,\mathbf{x},\mathbf{u},\psi_0,\boldsymbol{\psi}}
\newcommand{\xut}{t,\mathbf{x}(t),\mathbf{u}(t)}
\newcommand{\xu}{t,\mathbf{x},\mathbf{u}}
\newcommand{\maple}{\textsf{Maple} }

\newtheorem{theorem}{Theorem}
\newtheorem{example}[theorem]{Example}
\newtheorem{definition}[theorem]{Definition}
\newtheorem{remark}[theorem]{Remark}

\begin{document}

\title{Symbolic Computation of
       Variational Symmetries in Optimal Control\footnote{Presented at the
       4th Junior European Meeting on ``Control and Optimization'',
       Bia\l ystok Technical University, Bia\l ystok, Poland, 11-14 September 2005.
       Research report CM06/I-02. Accepted (24-Feb-2006) to Control \& Cybernetics.}}

\author{Paulo D. F. Gouveia\\
        \texttt{pgouveia@ipb.pt} \and
        Delfim F. M. Torres\\
        \texttt{delfim@mat.ua.pt} \and
        Eug\'{e}nio A. M. Rocha\\
        \texttt{eugenio@mat.ua.pt}}

\date{Control Theory Group (\textsf{cotg})\\
Centre for Research in Optimization and Control\\
Department of Mathematics, University of Aveiro\\
3810-193 Aveiro, Portugal}

\maketitle

\begin{abstract}
We use a computer algebra system to compute, in an efficient way,
optimal control variational symmetries up to a gauge term. The
symmetries are then used to obtain families of Noether's first
integrals, possibly in the presence of nonconservative external
forces. As an application, we obtain eight independent first integrals
for the sub-Riemannian nilpotent problem $(2,3,5,8)$.
\end{abstract}

\smallskip


\smallskip


\smallskip


\medskip


\section{Introduction}

The concept of variational symmetry entered into optimal control
in the seventies of the twentieth century \cite{Djukic73}.
Variational symmetries, which keep an optimal control problem
invariant, are described mathematically in terms of a group of
parameter transformations: two transformations performed one after
another may be replaced by one transformation of the same family;
there exists an identity transformation; to each transformation
there exists an inverse one. Variational symmetries are very
useful in optimal control, but unfortunately their study is not
easy, requiring lengthy and cumbersome calculations
\cite{Torres04}.

Recently there has been an interest in the application of Computer
Algebra Systems to the study of control systems, and collections
of symbolical tools are being developed to help on the analysis
and solution of complex problems. The first computer algebra
package for computing the variational symmetries in the calculus
of variations was given in \cite{gouv04}; then extended to the
more general setting of optimal control \cite{GouveiaTorresCMAM}.

In this work we provide a new Maple package for the automatic
computation of variational symmetries and respective Noether's
first integrals in the calculus of variations and optimal control.
The present package generalize the previous results in
\cite{GouveiaTorresCMAM} by introducing two new possibilities: (i)
invariance symmetries up to a gauge term \cite{delfimEJC}; (ii)
presence of nonconservative external forces \cite{gast05}.
Moreover, the efficiency in computing the variational symmetries
is largely improved when we compare the running times with the
ones in \cite{GouveiaTorresCMAM}. With the improvements in the
efficiency of the package, we are now able, for the first time
in the literature, to obtain eight independent first integrals
for the nilpotent problem $(2,3,5,8)$ of sub-Riemannian geometry.


\section{Nonconservative forces}
\label{sec:PMP}

Without loss of generality, we consider the optimal control
problem in Lagrange form: to minimize an integral functional
\begin{equation}
\label{eq:funcionalCO} I[\mathbf{x}(\cdot),\mathbf{u}(\cdot)] =
\int_{a}^{b} L(\xut) \,\dd t
\end{equation}
subject to a control system described by a system of ordinary
differential equations of the form
\begin{equation}
\label{eq:sistCont}
\dot{\mathbf{x}}(t)= \boldsymbol{\varphi}(\xut) \, ,
\end{equation}
together with appropriate boundary conditions, not relevant for
the present study (the results of the paper are valid for
arbitrary boundary conditions). The Lagrangian $L: \mathbb{R}
\times \mathbb{R}^n \times \mathbb{R}^m \rightarrow \mathbb{R}$
and the velocity vector $\boldsymbol{\varphi}\!:\mathbb{R}
\times\mathbb{R}^n\times\mathbb{R}^m \rightarrow\mathbb{R}^n$ are
assumed to be continuously differentiable functions with respect
to all their arguments. The controls $\mathbf{u}\!:[a,
b]\rightarrow\Omega\subseteq\mathbb{R}^m$ are piecewise continuous
functions taking values on an open set $\Omega$; the state
variables $\mathbf{x}\!:[a, b] \rightarrow\mathbb{R}^n$
continuously differentiable functions.

The resolution of optimal control problems usually goes by
identifying the Pontryagin extremals \cite{Pontryagin62}. In
presence of nonconservative external forces
$\mathbf{F} : \mathbb{R} \times \mathbb{R}^n \times
\mathbb{R}^m \rightarrow \mathbb{R}^n$
the Pontryagin Maximum Principle (PMP)
takes the following form \cite{gast05}.

\begin{theorem}[PMP under a nonconservative force $\mathbf{F}$]
\label{thm:pMaxPont} If
$\left(\mathbf{x}(\cdot),\mathbf{u}(\cdot)\right)$ is a solution
of the optimal control problem
(\ref{eq:funcionalCO})-(\ref{eq:sistCont}) under presence of a
nonconservative force
$\mathbf{F}\left(t,\mathbf{x},\mathbf{u}\right)$, then there
exists a non-vanishing pair
$\left(\psi_0,\boldsymbol{\psi}(\cdot)\right)$, where $\psi_0 \leq
0$ is a constant and $\boldsymbol{\psi}(\cdot)$ a $n$-vectorial
piecewise $C^1$-smooth function with domain $[a, b]$, in such a
way the quadruple
$(\mathbf{x}(\cdot),\mathbf{u}(\cdot),\psi_0,$
$\boldsymbol{\psi}(\cdot))$ satisfy the following conditions
almost everywhere in $[a, b]$:
\renewcommand{\theenumi}{\roman{enumi}}
\renewcommand{\labelenumi}{\emph{(\theenumi)}}
\begin{enumerate}
\item the nonconservative Hamiltonian system
\begin{equation}
\label{eq:sistCont3}
\begin{cases}
\dot{\mathbf{x}}(t)\T=
\frac{\partial H}{\partial\boldsymbol{\psi}}(\xupt)\, ,\\
\dot{\boldsymbol{\psi}}(t)\T = - \frac{\partial H}{\partial
\mathbf{x}}(\xupt) +\mathbf{F}(t,\mathbf{x}(t),\mathbf{u}(t))\T \, ;
\end{cases}
\end{equation}
\item the maximality condition
\begin{equation}
\label{eq:condMax}
H(\xupt)=
\begin{array}[t]{c}
\emph{max}\\
\mathbf{v}\hspace{-0.1cm}\in\hspace{-0.1cm}\Omega
\end{array}
H(t,\mathbf{x}(t),\mathbf{v},\psi_0,\boldsymbol{\psi}(t)) \, ;
\end{equation}
\end{enumerate}
where the Hamiltonian $H$ is defined by
\begin{equation}
\label{eq:hamilt}
H(\xup)= \psi_0
L(\xu)+
\boldsymbol{\psi}\T   \cdot
\boldsymbol{\varphi}(\xu) \, .
\end{equation}
\end{theorem}

\begin{remark}
The right-hand side of the equations of the nonconservative
Hamiltonian system (\ref{eq:sistCont3}) represent a row-vector.
First equation in (\ref{eq:sistCont3}) is nothing more than the
control system (\ref{eq:sistCont}); the second equation is known
as the \emph{nonconservative adjoint system}.
\end{remark}

\begin{definition}
\label{def:nonc:extr} A quadruple
$\left(\mathbf{x}(\cdot),\mathbf{u}(\cdot),\psi_0,
\boldsymbol{\psi}(\cdot)\right)$ satisfying
Theorem~\ref{thm:pMaxPont} is said to be a \emph{nonconservative
extremal}. A nonconservative extremal is said to be \emph{normal}
when $\psi_0\neq0$, \emph{abnormal} when $\psi_0=0$.
\end{definition}

\begin{remark}
Since we are assuming $\Omega$ to be an open set,
the maximality condition \eqref{eq:condMax} implies the \emph{stationary
condition}
\begin{equation}
\label{eq:stat:cond}
\frac{\partial H}{\partial \mathbf{u}}(\xupt)
= \mathbf{0} \, , \quad t \in [a,b] \, .
\end{equation}
\end{remark}


\section{Invariance up to a gauge term}
\label{sec:gauge}

Let $\mathbf{h}^s:[a, b] \times \mathbb{R}^n \times \mathbb{R}^m
\times \mathbb{R} \times \mathbb{R}^n \rightarrow \mathbb{R}
\times \mathbb{R}^n \times \mathbb{R}^m \times \mathbb{R}^n$ be a
one-parameter group of $\mathbb{C}^1$ transformations of the form
\begin{eqnarray}
\label{eq:transf}
\mathbf{h}^s(\xup)= \hspace{8.1cm} \nonumber\\
(h_t^s(\xup),\mathbf{h}_\mathbf{x}^s(\xup),
\mathbf{h}_\mathbf{u}^s(\xup),
\mathbf{h}_{\boldsymbol{\psi}}^s(\xup)) \, .
\end{eqnarray}
Without loss of generality, we assume that the identity
transformation of the group (\ref{eq:transf}) is obtained when the
parameter $s$ is zero:
\begin{eqnarray*}
h_t^0(\xup)=t,\;
\mathbf{h}_\mathbf{x}^0(\xup)=
\mathbf{x},\;\nonumber\\
\mathbf{h}_\mathbf{u}^0(\xup)=
\mathbf{u},\;
\mathbf{h}_{\boldsymbol{\psi}}^0(\xup)=\boldsymbol{\psi}.
\end{eqnarray*}
Associated with a one-parameter group of transformations
(\ref{eq:transf}), we introduce its \emph{infinitesimal generators}:
\begin{eqnarray}
\label{eq:transf2}
T(\xup)  &=&  \left. \frac{\partial}
{\partial{s}} h^s_t\right|_{s=0}\textrm{, }
\mathbf{X}(\xup)  = \left.
\frac{\partial}{\partial{s}} \mathbf{h}_\mathbf{x}^s\right|_{s=0}
\textrm{,}\nonumber\\
\mathbf{U}(\xup)  &=&
\left. \frac{\partial} {\partial{s}} \mathbf{h}_\mathbf{u}^s
\right|_{s=0} \textrm{, }
\boldsymbol{\Psi}(\xup)  =  \left.
\frac{\partial}{\partial{s}}
\mathbf{h}_{\boldsymbol{\psi}}^s\right|_{s=0} \textrm{. }
\end{eqnarray}

\begin{definition}[Invariance up to a gauge term]
\label{defin:invarGauge} An optimal control problem
(\ref{eq:funcionalCO})-(\ref{eq:sistCont}) is said to be invariant
under a one-parameter group of transformations (\ref{eq:transf})
up to a gauge term $g^s(t,\mathbf{x},\mathbf{u},\psi_0,
\boldsymbol{\psi}) \in \mathbb{C}^1([a, b], \mathbb{R}^n,
\mathbb{R}^m, \mathbb{R}, \mathbb{R}^n; \mathbb{R})$, if for all
$s$ sufficiently small and for any subinterval $[\alpha,\beta]
\subseteq [a,b]$ one has
\begin{eqnarray}
\label{eq:invar}
\int_{\alpha^s}^{\beta^s}\left(H(t^s,\mathbf{x}^s(t^s),\mathbf{u}^s(t^s),
\psi_0,\boldsymbol{\psi}^s(t^s))
-\boldsymbol{\psi}^s(t^s)\T \cdot \frac{\dd}{\dd t^s}\mathbf{x}^s(t^s) \right)
\,\emph{d}t^s \nonumber\\
=\int_{\alpha}^{\beta}\biggl(H(\xupt)
-\boldsymbol{\psi}(t)\T  \cdot  \frac{\dd}{\dd t}\mathbf{x}(t) \hspace{2.5cm}\nonumber\\
+\frac{\dd}{\dd t} g^s(\xupt)\biggr)\,\dd t\, ,
\end{eqnarray}
where
$\alpha^s=h_t^s(\alpha,\mathbf{x}(\alpha),\mathbf{u}(\alpha),\psi_0,
\boldsymbol{\psi}(\alpha))$,
$\beta^s=h_t^s(\beta,\mathbf{x}(\beta),\mathbf{u}(\beta),\psi_0,
\boldsymbol{\psi}(\beta))$, and
$\left(t^s,\mathbf{x}^s,\mathbf{u}^s,\boldsymbol{\psi}^s\right) =
\left(h_t^s,\mathbf{h}_\mathbf{x}^s,\mathbf{h}_\mathbf{u}^s,
\mathbf{h}_{\boldsymbol{\psi}}^s\right)$.
\end{definition}

When we write \eqref{eq:invar} in terms of the generators
\eqref{eq:transf2}, one gets a necessary and sufficient condition
of invariance -- \textrm{cf.} \cite{Djukic73,Torres05}.

\begin{theorem}[Necessary and sufficient condition of invariance]
\label{thm:condInv} An optimal control problem is \emph{invariant}
under \eqref{eq:transf2} up to $G(\xup) =\left. \frac{\dd}{\dd s}
g^s(t,\mathbf{x},\mathbf{u},\psi_0,
\boldsymbol{\psi})\right|_{s=0}$ or, equivalently,
\eqref{eq:transf2} is a \emph{symmetry} of the problem up to $G$,
if, and only if,
\begin{equation}
\label{eq:detGerad}
\frac{\partial H}{\partial t}T +\frac{\partial
H}{\partial \mathbf{x}}\cdot \mathbf{X} +\frac{\partial
H}{\partial \mathbf{u}}\cdot \mathbf{U} +\frac{\partial
H}{\partial \boldsymbol{\psi}}\cdot \boldsymbol{\Psi}
-\boldsymbol{\Psi}\T \cdot \dot{\mathbf{x}}
-\boldsymbol{\psi}\T \cdot \frac{\dd
\mathbf{X}}{\dd t} +H \frac{\dd T}{\dd t}=
\frac{\dd G}{\dd t}\, ,
\end{equation}
with $H$ the Hamiltonian \eqref{eq:hamilt}.
\end{theorem}
\begin{remark}
The function $G(\xup) =\left. \frac{\dd}{\dd s}
g^s(t,\mathbf{x},\mathbf{u},\psi_0,
\boldsymbol{\psi})\right|_{s=0}$ is also known in the literature
as a \emph{gauge term}.
\end{remark}
\begin{proof}
Transforming the integral on the left-hand side of
(\ref{eq:invar}) to the interval $[\alpha, \beta]$, and having in
mind that (\ref{eq:invar}) is satisfied for all subintervals
$[\alpha,\beta] \subseteq [a,b]$, the invariance condition can be
written in the following equivalent form:
\begin{eqnarray*}
\left(H(\mathbf{h}^s(t,\mathbf{x},
\mathbf{u},\psi_0,\boldsymbol{\psi}))
-\mathbf{h}_{\boldsymbol{\psi}}^s(t,\mathbf{x},\mathbf{u},\psi_0,
\boldsymbol{\psi})\T \cdot \frac{\frac{\dd\mathbf{h}_\mathbf{x}^s(\xup)}
{\dd t}}{\frac{\dd h_t^s(\xup)}{\dd t}
}
 \right)
\frac{\textrm{d}h_t^s(t,\mathbf{x},\mathbf{u},\psi_0, \boldsymbol{\psi})}
{\dd t} \\
= H(\xup) -\boldsymbol{\psi} \T \cdot \frac{\dd}{\dd t}\mathbf{x}
+ \frac{\dd}{\dd t} g^s(\xup) \, .
\end{eqnarray*}
Differentiating  both sides of the equation with respect to $s$,
\begin{eqnarray*}
\frac{\dd}{\dd s} \Biggl[\left( H(\mathbf{h}^s(\xup))
-\mathbf{h}_{\boldsymbol{\psi}}^s(\xup)\T \cdot
\frac{\frac{\dd \mathbf{h}_\mathbf{x}^s(\xup)}{\dd t}
}
{\frac{\dd h_t^s(\xup)}{\dd t}
}
 \right)\\
 \times \frac{\dd h_t^s(\xup)}{\dd t}
 \Biggr]
= \frac{\dd}{\dd s} \left(\frac{\dd}{\dd t} g^s(\xup)\right),
\end{eqnarray*}
we obtain the equality
\begin{eqnarray*}
\left( H(\mathbf{h}^s)
-\mathbf{h}_{\boldsymbol{\psi}}^s{}\T \cdot
\frac{\dd \mathbf{h}_\mathbf{x}^s/\dd t}
{\dd h_t^s/\dd t}
 \right)
 \frac{\dd}{\dd t}\frac{\dd h_t^s}{\dd s}
 +\Biggl(
 \frac{\partial H(\mathbf{h}^s)}{\partial h_t^s}\frac{\partial h_t^s}
 {\partial s}
 + \frac{\partial H(\mathbf{h}^s)}{\partial \mathbf{h}_\mathbf{x}^s}
 \cdot \frac{\partial \mathbf{h}_\mathbf{x}^s}{\partial s}
 \\
 + \frac{\partial H(\mathbf{h}^s)}{\partial \mathbf{h}_\mathbf{u}^s}
 \cdot \frac{\partial \mathbf{h}_\mathbf{u}^s}{\partial s}
 + \frac{\partial H(\mathbf{h}^s)}{\partial \mathbf{h}_{\boldsymbol{\psi}}^s}
 \cdot \frac{\partial \mathbf{h}_{\boldsymbol{\psi}}^s}{\partial s}
-\frac{\dd {\mathbf{h}_{\boldsymbol{\psi}}^s}\T}{\dd s} \cdot
\frac{\dd \mathbf{h}_\mathbf{x}^s/\dd t}
{\dd h_t^s/\dd t}
\\
-\mathbf{h}_{\boldsymbol{\psi}}^s{}\T \cdot
\left(
\frac{\frac{\dd}{\dd t}\frac{\dd\mathbf{h}_\mathbf{x}^s}{\dd s}}
{\frac{\dd h_t^s}{\dd t}}
-\frac{\frac{\dd\mathbf{h}_\mathbf{x}^s}{\dd t}\frac{\dd}{\dd t}
\frac{\dd h_t^s}{\dd s}} {\frac{\dd h_t^s}{\dd t}\frac{\dd h_t^s}{\dd t}}
\right)
 \Biggr)
 \frac{\dd h_t^s}{\dd t}
 =
 \frac{\dd}{\dd t} \frac{\dd g^s}{\dd s} \, .
\end{eqnarray*}
Finally, choosing $s=0$, we express the condition in terms of the
infinitesimal generators \eqref{eq:transf2} and the function
$G(\xup) =\left. \frac{\dd}{\dd s}
g^s(t,\mathbf{x},\mathbf{u},\psi_0,
\boldsymbol{\psi})\right|_{s=0}$:
\begin{eqnarray*}
\left( H
-\boldsymbol{\psi}{}\T \cdot
\dot{\mathbf{x}} \right)
 \frac{\dd T}{\dd t}
+\biggl(
 \frac{\partial H}{\partial t} T
 + \frac{\partial H}{\partial \mathbf{x}}
 \cdot \mathbf{X}
 + \frac{\partial H}{\partial \mathbf{u}}
 \cdot \mathbf{U}
 + \frac{\partial H}{\partial \boldsymbol{\psi}}
 \cdot \boldsymbol{\Psi}
-\boldsymbol{\Psi}\T \cdot
\dot{\mathbf{x}}\\
-\boldsymbol{\psi}\T \cdot
\left(
\frac{\dd\mathbf{X}}{\dd t}
-\dot{\mathbf{x}} \frac{\dd T}{\dd t}
\right)
 \biggr)
=
 \frac{\dd G}{\dd t} \, .
\end{eqnarray*}
\end{proof}


\section{Nonconservative Noether's theorem}
\label{sec:noncons}

Emmy Noether was the first who established the relation between
the existence of invariance transformations of the problems and
the existence of conservation laws -- first integrals of the
Euler-Lagrange or Hamiltonian equations \cite{Noether}. A
generalization of the classical result of E.~Noether for the
nonconservative calculus of variations was recently given by Fu
and Chen \cite{Fu03}; then extended to the more general setting of
optimal control by Frederico and Torres \cite{gast05}.

Using (\ref{eq:sistCont3}), together with the stationary condition
(\ref{eq:stat:cond}), one can deduce that along the
nonconservative Pontryagin extremals
(Definition~\ref{def:nonc:extr}), the total derivative of the
Hamiltonian with respect to the independent variable $t$ is equal
to its partial derivative plus the scalar product of the velocity
vector with the resultant nonconservative forces $\mathbf{F}$
\cite{gast05}:
\begin{equation}
\label{eq:prop:dHdt2} \frac{\dd}{\dd t} H(\xupt) =
\frac{\partial}{\partial t}H(\xupt) + \dot{\mathbf{x}}(t)\T \cdot
\mathbf{F}(t,\mathbf{x}(t),\mathbf{u}(t)) \, .
\end{equation}
Using this fact, the nonconservative optimal control version of
E.~Noether's theorem is easily obtained from the necessary and
sufficient invariance condition (\ref{eq:detGerad}), restricting
attention to the quadruples
$(\mathbf{x}(\cdot),\mathbf{u}(\cdot),\psi_0,
\boldsymbol{\psi}(\cdot))$ that satisfy the nonconservative
Hamiltonian system (\ref{eq:sistCont3}) and the maximality
condition (\ref{eq:condMax}): along the extremals, equalities
(\ref{eq:sistCont3}), \eqref{eq:stat:cond}, and
\eqref{eq:prop:dHdt2} permit to simplify \eqref{eq:detGerad} to
the form
\begin{eqnarray*}
&&\left(\frac{\dd H}{\dd t} - \dot{\mathbf{x}}\T
\cdot \mathbf{F} \right)T
+\left(\mathbf{F}\T
- \dot{\boldsymbol{\psi}}\T\right) \cdot \mathbf{X} -
\boldsymbol{\psi}\T \cdot \frac{\dd
\mathbf{X}}{\dd t} + H \frac{\dd T}{\dd t}=
\frac{\textrm{d}G}{\textrm{d}t}\nonumber\\
&& \Leftrightarrow
\frac{\dd H}{\dd t} T + H \frac{\dd T}{\dd t}
- \dot{\boldsymbol{\psi}}\T \cdot \mathbf{X} -
\boldsymbol{\psi}\T \cdot \frac{\dd
\mathbf{X}}{\dd t}
-\frac{\textrm{d}G}{\textrm{d}t}
- \left(\dot{\mathbf{x}}\T\,T
-\mathbf{X}\T\right) \cdot \mathbf{F}=0\nonumber\\
&& \Leftrightarrow
\frac{\dd}{\dd t} \left( H T -
\boldsymbol{\psi}\T \cdot \mathbf{X} - G
- \int\!\!\left(\dot{\mathbf{x}}\T\,T
-\mathbf{X}\T\right) \cdot \mathbf{F}\,\dd t
\right) = 0 \, .
\end{eqnarray*}
This means that $H T - \boldsymbol{\psi}\T \cdot \mathbf{X} - G -
\int\!\!\left(\dot{\mathbf{x}}\T\,T -\mathbf{X}\T\right) \cdot
\mathbf{F}\,\dd t$ is a first integral whenever the optimal
control problem under consideration admits a symmetry
\eqref{eq:transf2} up to the gauge term $G$:

\begin{theorem}[Nonconservative Optimal Control version of Noether's Principle]
\label{thm:TNoether2} If the infinitesimal generators
(\ref{eq:transf2}) constitute a symmetry of the optimal control
problem (\ref{eq:funcionalCO})-(\ref{eq:sistCont}) under presence
of nonconservative forces with resultant vector
$\mathbf{F}\left(t,\mathbf{x},\mathbf{u}\right)$, then
\begin{eqnarray}
\label{eq:leinaocons}
\int\!\!\left(\dot{\mathbf{x}}(t)\T T(\xupt)
\!-\!\!\mathbf{X}(\xupt)\T\right)
\!\cdot\! \mathbf{F}(t,\mathbf{x}(t),\mathbf{u}(t))
\dd t\nonumber\\
+\;\boldsymbol{\psi}(t)\T \cdot \mathbf{X}(\xupt)+G(\xupt)
\nonumber\\
-H(\xupt)\;T(\xupt) = \textrm{const}
\nonumber\\
\end{eqnarray}
is a conservation law, \textrm{i.e.}, condition
\eqref{eq:leinaocons} holds for all $t$ in $[a,b]$ and for every
nonconservative extremal
$\left(\mathbf{x}(\cdot),\mathbf{u}(\cdot),\psi_0,
\boldsymbol{\psi}(\cdot)\right)$ of the problem.
\end{theorem}


\section{Computation of symmetries up to a gauge term}
\label{sec:comput}

The main problem in obtaining Noether's conservation laws (in
applying Theorem~\ref{thm:TNoether2}) resides in the determination
of the symmetries and respective gauge terms. If $n$ effective
first integrals exist \cite{RochaTorres4thJM}, then the optimal
control problem is integrable, and classical results allow the
integration of the equations of motion.

Here we propose an algorithm for determining the infinitesimal
generators \eqref{eq:transf2} and the gauge terms $G$ which define
a variational symmetry. Let us assume, for the moment, that the
optimal controls are $C^1$ functions (in \S\ref{sec:optdep} we
will drop this restrictive assumption, just by assuming that $T$,
$\mathbf{X}$, and $G$ do not depend on the control variables). The
key point to compute symmetries consists in generalizing the method used
in \cite[\S3]{GouveiaTorresCMAM} to the nonconservative and
gauge-invariant cases. The idea is simple: when we substitute the
Hamiltonian $H$ and its partial derivatives in the invariance
identity \eqref{eq:detGerad}, then the condition becomes a
polynomial in $\dot{\mathbf{x}}$, $\dot{\mathbf{u}}$ and
$\dot{\boldsymbol{\psi}}$, and one can equal the coefficients of
the polynomial to zero. Thus, given an optimal control problem
(\ref{eq:funcionalCO})-(\ref{eq:sistCont}), defined by a
Lagrangian $L$ and a velocity vector $\boldsymbol{\varphi}$, we
determine the infinitesimal generators $T$, $\mathbf{X}$,
$\mathbf{U}$ and $\boldsymbol{\Psi}$ and the gauge term $G$, which
define a symmetry for the problem, by the following method: (i) we
define the respective Hamiltonian \eqref{eq:hamilt}; (ii) we
substitute $H$ and its partial derivatives into
\eqref{eq:detGerad}; (iii) expanding the total derivatives
\begin{eqnarray}
\frac{\textrm{d}T}{\textrm{d}t} &=& \frac{\partial T}{\partial t}
+\frac{\partial T}{\partial \mathbf{x}}\cdot\dot{\mathbf{x}}
+\frac{\partial T}{\partial \mathbf{u}}\cdot\dot{\mathbf{u}}
+\frac{\partial T}{\partial
\boldsymbol{\psi}}\cdot\dot{\boldsymbol{\psi}}\textrm{,}\quad \notag \\
\frac{\textrm{d}\mathbf{X}}{\textrm{d}t} &=& \frac{\partial
\mathbf{X}}{\partial t}
+ \frac{\partial \mathbf{X}} {\partial\mathbf{x}}\cdot\dot{\mathbf{x}}
+ \frac{\partial \mathbf{X}} {\partial\mathbf{u}}\cdot\dot{\mathbf{u}}
+ \frac{\partial \mathbf{X}}{\partial \boldsymbol{\psi}}\cdot
\dot{\boldsymbol{\psi}}, \label{eq:smothCtrl} \\
\frac{\textrm{d}G}{\textrm{d}t} &=& \frac{\partial G}{\partial t}
+\frac{\partial G}{\partial \mathbf{x}}\cdot\dot{\mathbf{x}}
+\frac{\partial G}{\partial \mathbf{u}}\cdot\dot{\mathbf{u}}
+\frac{\partial G}{\partial
\boldsymbol{\psi}}\cdot\dot{\boldsymbol{\psi}}\textrm{,} \notag
\end{eqnarray}
we write equation (\ref{eq:detGerad}) as a polynomial
\begin{equation}
\label{eq:poly}
A(\xup)+B(\xup) \cdot \dot{\mathbf{x}} +
C(\xup) \cdot \dot{\mathbf{u}} + D(\xup)  \cdot \dot{\boldsymbol{\psi}} = 0
\end{equation}
in the $2n+m$ derivatives $\dot{\mathbf{x}}$, $\dot{\mathbf{u}}$ and
$\dot{\boldsymbol{\psi}}$:
\begin{eqnarray}
\label{eq:detGerad2}
\left(\frac{\partial H}{\partial t} T
+\frac{\partial H}{\partial \mathbf{x}}\cdot \mathbf{X}
+\frac{\partial H}{\partial \mathbf{u}}\cdot \mathbf{U}
+\frac{\partial H}{\partial \boldsymbol{\psi}}\cdot \boldsymbol{\Psi}
+H \frac{\partial T}{\partial t}
-\boldsymbol{\psi}\T \cdot
\frac{\partial \mathbf{X}}{\partial t}
- \frac{\partial G}{\partial t}
\right)\nonumber\\
+ \left(
-\boldsymbol{\Psi}\T
+H \frac{\partial T}{\partial\mathbf{x}}
-\boldsymbol{\psi}\T \cdot \frac{\partial \mathbf{X}}
{\partial \mathbf{x}}
-\frac{\partial G}{\partial \mathbf{x}}
\right)\cdot \dot{\mathbf{x}}
+\left(
H \frac{\partial T}{\partial \mathbf{u}}
-\boldsymbol{\psi}\T \cdot
\frac{\partial \mathbf{X}}{\partial \mathbf{u}}
-\frac{\partial G}{\partial \mathbf{u}}
\right) \cdot \dot{\mathbf{u}}
\nonumber\\
+\left(
H \frac{\partial T}{\partial \boldsymbol{\psi}}
-\boldsymbol{\psi}\T \cdot
\frac{\partial \mathbf{X}}{\partial \boldsymbol{\psi}}
-\frac{\partial G}{\partial\boldsymbol{\psi}}
\right)
\cdot \dot{\boldsymbol{\psi}}=0\, . \quad
\end{eqnarray}
The terms in \eqref{eq:detGerad2}, which involve derivatives with
respect to vectors, are expanded in row-vectors or in matrices,
depending, respectively, if the function is a scalar or a
vectorial one. For example,
\begin{eqnarray}
\frac{\partial T}{\partial \mathbf{x}}
&=&\left[ \frac{\partial T}{\partial x_1}
\ \frac{\partial T}{\partial x_2}
\ \cdots \ \frac{\partial T}{\partial x_n}\right],
\nonumber\\
\frac{\partial \mathbf{X}}{\partial \boldsymbol{\psi}}&=&
\left[\frac{\partial \mathbf{X}}{\partial \psi_1}\ \frac{\partial
\mathbf{X}}{\partial \psi_2}\ \cdots\  \frac{\partial
\mathbf{X}}{\partial \psi_n}\right] = \left[
\begin{array}{cccc}
\smallskip
\frac{\partial X_1}{\partial \psi_1} & \frac{\partial X_1}{\partial \psi_2}
& \cdots & \frac{\partial X_1}{\partial \psi_n}\\
\frac{\partial X_2}{\partial \psi_1} & \frac{\partial X_2}{\partial \psi_2}
& \cdots & \frac{\partial X_2}{\partial \psi_n}\\
\vdots&\vdots&\ddots&\vdots\\
\frac{\partial X_n}{\partial \psi_1} & \frac{\partial X_n}{\partial \psi_2}
& \cdots & \frac{\partial X_n}{\partial \psi_n}
\end{array}
\right].\nonumber
\end{eqnarray}
 Equation
(\ref{eq:detGerad2}) is a differential equation in the $2n+m+2$
unknown functions $T$, $X_1$, \dots, $X_n$, $U_1$, \dots, $U_m$,
$\varPsi_1$, \dots, $\varPsi_n$ and $G$. This equation must hold
for all $\dot{x}_1$, \dots, $\dot{x}_n$, $\dot{u}_1$, \dots,
$\dot{u}_n$, $\dot{\psi}_1$, \dots, $\dot{\psi}_n$, and therefore
the coefficients $A$, $B$, $C$ and $D$ of polynomial
\eqref{eq:poly} must vanish, that is,
\begin{equation}
\label{eq:detGerad3}
\left\{
\begin{array}{l}
\displaystyle \medskip
\frac{\partial H}{\partial t}T
+\frac{\partial H}{\partial \mathbf{x}}\cdot \mathbf{X}
+\frac{\partial H}{\partial \mathbf{u}}\cdot \mathbf{U}
+\frac{\partial H}{\partial \boldsymbol{\psi}}\cdot \boldsymbol{\Psi}
+H \frac{\partial T}{\partial t}
-\boldsymbol{\psi}\T \cdot
\frac{\partial \mathbf{X}}{\partial t}- \frac{\partial G}{\partial t}=0 \, ,\\
\displaystyle \medskip
-\boldsymbol{\Psi}\T
+H \frac{\partial T}{\partial\mathbf{x}}
-\boldsymbol{\psi}\T \cdot \frac{\partial \mathbf{X}}
{\partial \mathbf{x}}
-\frac{\partial G}{\partial \mathbf{x}}
=\mathbf{0} \, , \\
\displaystyle \medskip
H \frac{\partial T}{\partial \mathbf{u}}
-\boldsymbol{\psi}\T \cdot
\frac{\partial \mathbf{X}}{\partial \mathbf{u}}
-\frac{\partial G}{\partial \mathbf{u}}
=\mathbf{0} \, ,\\
\displaystyle
H \frac{\partial T}{\partial \boldsymbol{\psi}}
-\boldsymbol{\psi}\T \cdot
\frac{\partial \mathbf{X}}{\partial \boldsymbol{\psi}}
-\frac{\partial G}{\partial\boldsymbol{\psi}}
=\mathbf{0} \, .
\end{array}
\right.
\end{equation}
System of equations \eqref{eq:detGerad3}, obtained from
(\ref{eq:detGerad2}), is a system of $2n+m+1$ partial differential
equations with $2n+m+2$ unknown functions so, in general, there
exists not a unique symmetry but a family of such symmetries.
The system \eqref{eq:detGerad3} becomes even more under-determined when one assumes,
as in \S\ref{sec:optdep}, that $T$, $\mathbf{X}$, and $G$ do not depend on the
control variables $\mathbf{u}$.
Although a system of partial differential equations,
solving \eqref{eq:detGerad3} is possible because
the system is of the first order
and linear with respect to the unknown functions and their
derivatives. We solve the system of PDEs by the method
of (additive) separation of variables, as explained in
\cite{Terrab95}. Following \cite{Terrab95}, the generators are
replaced by the sum of unknown functions, one for each variable.
For example, $T(t,x_1,x_2,\psi_1,\psi_2) =
T_1(t)+T_2(x_1)+T_3(x_2)+T_4(\psi_1)+T_5(\psi_2)$.
When dealing with optimal control problems
with several state and control variables, the number of
calculations is big enough, and the help of the computer is more
than welcome. We define a \maple procedure \emph{Symmetry} that
does all the cumbersome calculations for us.
The procedure receives, as input, the Lagrangian and the velocity
vector; and returns, as output, a family of symmetries
$\left(T,\mathbf{X},\mathbf{U},\boldsymbol{\Psi}\right)$ and, if
necessary, the respective gauge term $G$.
We remark that since system \eqref{eq:detGerad3} is homogeneous,
we always have, as trivial solution,
$\left(T,\mathbf{X},\mathbf{U},\mathbf{\Psi}\right) = \mathbf{0}$.


\section{The computer algebra package}
\label{sec:update}

We obtain Noether conservation laws, in an automatic way, through
two steps: (i) with our procedure \emph{Symmetry} we obtain the
variational symmetries and respective gauge terms; (ii) using the
obtained symmetries, gauge terms, and nonconservative forces as
input to procedure \emph{Noether}, we obtain the correspondent
conservation laws. In \S\ref{sec:exOC} we give several examples,
not covered by the previous results in
\cite{gouv04,GouveiaTorresCMAM}, illustrating the whole process.
Given the limit on the maximum number of pages of the paper, we do
not provide the \maple definitions for the procedures
\emph{Symmetry} and \emph{Noether} here. The complete \maple
package can be freely obtained from
\texttt{http://www.mat.ua.pt/delfim/maple.htm}\, together with an
online help database for the \maple system.

Novelties of the procedures \emph{Symmetry} and \emph{Noether}
with respect to the previous versions in
\cite{gouv04,GouveiaTorresCMAM} are: (i) possibility of procedure
\emph{Symmetry} to cover invariance symmetries up to a gauge term,
according with \S\ref{sec:gauge} and \S\ref{sec:comput}; (ii)
improvements of efficiency -- see \S\ref{sec:optdep}; (iii)
possibility of procedure \emph{Noether} to consider problems of
the calculus of variations and optimal control under
nonconservative external forces, according with
\S\ref{sec:noncons}; (iv) improvement of the usage of the
procedures by introduction of several optional parameters, as
illustrated in \S\ref{sec:exOC}. Moreover, a new \maple procedure
called \emph{PMP} was added which implements
Theorem~\ref{thm:pMaxPont}, according with \S\ref{sec:PMP}.\footnote{In
the software Cotcot, available from
\textsf{http://www.n7.fr/apo/cotcot/}, the tool Adifor
for automatic differentiation of Fortran is also used
to generate, in the conservative case,
the equations of the Pontryagin maximum principle \cite{cotcot}.} The
procedure \emph{PMP} is very useful in practice, when dealing with
concrete problems of the calculus of variations and optimal
control -- \textrm{cf.} \S\ref{sec:exOC}. The input to the
procedure is: the Lagrangian $L$ and the velocity vector
$\boldsymbol{\varphi}$, that define the optimal control problem
\eqref{eq:funcionalCO}-\eqref{eq:sistCont} and the respective
Hamiltonian $H$; the nonconservative external forces (if present);
and several useful optional arguments which define the output. The
output of \emph{PMP} is either (depending on the optional
parameters): the (nonconservative) extremals; the equations of the
(nonconservative) Hamiltonian system and stationary condition; or,
alternatively, the Hamiltonian. We refer the reader to the
Examples on \S\ref{sec:exOC} for a general overview on the usage
of the developed \maple procedures; to the annotated \maple
worksheet available at
\textsf{http://www.mat.ua.pt/delfim/maple.htm}, with all the
definitions of the package, detailed documentation, and many other
examples not given here, for more details. The reader is free to
experiment the \maple package in order to determine variational
symmetries and Noether conservation laws on his/her own problems.


\section{Efficiency, comparison with previous results}
\label{sec:optdep}

The high number of dependences that the infinitesimal generators
may present, affect, excessively, the efficiency of the method
described in \S\ref{sec:comput}, namely for problems with a large
number of state and control variables. In order to quantify this
effect, we measured the computing running times of our procedure
\emph{Symmetry} for different dependences of the infinitesimal
generators \eqref{eq:transf2}, with a large set of optimal control
problems: the ten problems considered in \cite[\S4,
\S5]{GouveiaTorresCMAM} (examples 4.1--4.6 and 5.1--5.4), together
with twelve new problems. Three of these new problems are given in
\S\ref{sec:exOC}, the complete set of problems being available as
a \maple worksheet, as mentioned in \S\ref{sec:update}. All the
computational processing was carried out with the \maple10
Computer Algebra System on a 1.4GHz Pentium Centrino with 512MB of
RAM. In the previous work \cite{GouveiaTorresCMAM}, the maximum
number of dependences for each generator, as indicated in
(\ref{eq:transf2}), is always considered. We denote here such
situation by \emph{D1}. In the \emph{D1} case, and as noticed in
\cite{GouveiaTorresCMAM}, the involved computational effort is
sometimes very high: the computing times increase exponentially
with the dimension of the problem. This is particularly well
illustrated with the following problems of sub-Riemannian
geometry: nilpotent problem $(2,3)$, with three state variables,
requires a total computing time of one minute (\cite[Example
4.5]{GouveiaTorresCMAM}); problem $(2, 3, 5)$, with five state
variables, requires thirty minutes
(\cite[Example~4.6]{GouveiaTorresCMAM}); the problem $(2,3,5,8)$,
with eight state variables, was not studied in
\cite{GouveiaTorresCMAM}, and thought to be out of its capacities.
We compute here its symmetries in Example~\ref{ex:prob2358}, with
the present \maple package, with forty one minutes of computing
time; while the method in \cite{GouveiaTorresCMAM} requires,
approximately, thirty times this value: twenty hours of computing
time are needed.\footnote{We believe that the forty minutes of computing
time can still be diminished by using a programming language closer to machine,
for instance using Adifor: \textsf{http://www-unix.mcs.anl.gov/autodiff/ADIFOR}.}

The computing running times largely depend on the numbers $n$ and
$m$, respectively the number of state and control variables:
besides directly influencing the number of dependences of the
unknown functions (infinitesimal generators), they determine the
amount of those functions and the number of partial differential
equations that must be solved in order to find the variational
symmetries. Without considering the gauge term, we come across a
system of $m+2n+1$ partial differential equations and $m+2n+1$
unknown functions, each one of the unknown functions being
dependent of $m+2n+1$ variables. We address here the following
question: is there some way to simplify the process of obtaining
the variational symmetries?

Although knowing that the complexity of the method is intimately
related with the values $n$ and $m$, that are fixed with a given
optimal control problem, we get, even so, a quite satisfactory
answer to the question. Analyzing the results from the test set of
problems, we verify that, in spite of considering the maximum
number of dependences (\emph{D1}), the infinitesimal generators
obtained through the procedure \emph{Symmetry} are, nevertheless,
almost always, dependent functions of a quite reduced number of
variables. When we restrict ourselves to the dependences $T(t)$,
$\mathbf{X}(t,\mathbf{x})$,
$\mathbf{U}(\mathbf{u},\boldsymbol{\psi})$,
$\boldsymbol{\Psi}(\boldsymbol{\psi})$ -- that we identify as
\emph{D2} -- we are able to cover the totality of the twenty two
considered problems in our study. If in the formulation of the
system of PDEs \eqref{eq:detGerad3} we only enter with these
dependences, besides the obvious reduction of the number of
dependences of the unknown functions, we reduce the number of
equations for less of half: from $m+2n+1$ to $n+1$. In agreement
with the simulations done, the efficiency of the procedure
\emph{Symmetry} increases significantly with this new group of
dependences (\emph{D2}). For instance, for the problem $(2,3,5)$
of sub-Riemannian geometry \cite[Example~4.6]{GouveiaTorresCMAM},
a problem with two controls and five state variables, the running
time passed from half an hour to less than one and a half minute.
We have also considered another more simplified set of
dependences, denoted by \emph{D3}: $T(t)$,
$\mathbf{X}(t,\mathbf{x})$, $\mathbf{U}(t,\mathbf{u})$,
$\boldsymbol{\Psi}(t,\boldsymbol{\psi})$. With it, is now possible
to obtain the symmetries of the sub-Riemannian nilpotent problem
$(2,3,5,8)$ (Example~\ref{ex:prob2358}), in less than 45 minutes;
and it is still possible to obtain the same conservation laws for
all the twenty two studied problems (in three of the problems
\cite[Examples~4.4, 5.2 and 5.3]{GouveiaTorresCMAM} the generators
were different, since the more general generators $\mathbf{U}$
depend on the variables $\boldsymbol{\psi}$, but the correspondent
Noether conservation laws \eqref{eq:leinaocons} are exactly the
same since they only depend on the generators $T$ and
$\boldsymbol{X}$). Finally, we repeated the study for a more
restricted group of dependences (\emph{D4}): $T(t)$,
$\mathbf{X}(\mathbf{x})$, $\mathbf{U}(\mathbf{u})$,
$\boldsymbol{\Psi}(\boldsymbol{\psi})$. As expected, the time of
processing suffered an additional reduction (for the $(2,3,5,8)$
problem the running time passed from $44'16 ''$ to $28'21 ''$),
but in this case not all the family of conservation laws for the
problems are obtained. For four of the problems --
\cite[Example~4.3]{GouveiaTorresCMAM}, Examples~2 and 3 in the \maple worksheet,
and Example~\ref{ex:Ex1_1} -- only particular cases of the
complete family of conservation laws are obtained.

To summarize the influence that the different dependences of the
generators have in the efficiency of the procedure
\emph{Symmetry}, we give in Table~\ref{tab:depGerad} the running
times for computing the variational symmetries of the three
problems of sub-Riemannian geometry already mentioned:
\cite[Examples~4.5 and 4.6]{GouveiaTorresCMAM} and
Example~\ref{ex:prob2358}. All the three problems have two control
variables and the same Lagrangian, but a different number of state
variables, respectively, 3, 5, and 8.

\begin{table}[!h]
\begin{center}
\small
\begin{tabular}{|c|r|r|r|r|}
\hline
{Dependences} & Nº PDEs${}^{*}$& prob. (2, 3) &
prob. (2, 3, 5) & prob. (2, 3, 5, 8)\\
\hline \hline
\emph{D1}& $m\!+\!2n\!+\!1$ & $1'04''$ & $30'34''$ &  $20h07'12''$\\
\hline
\emph{D2}& $n+1$ & $5''$  & $1'26''$ & $51'28''$\\
\hline
\emph{D3}& $n+1$ & $4''$ &  $1'09''$ & $44'16''$\\
\hline
\emph{D4}& $n+1$ & $2''$ & $38''$ & $28'21''$\\
\hline \multicolumn{5}{p{9cm}}{${}^{*}$ $n=$ nº of state
variables;
$m=$ nº of control variables.} \\
\end{tabular}
\normalsize
\caption{Running times of procedure \emph{Symmetry}
for three problems of sub-Riemannian geometry (\cite[Examples~4.5,
4.6]{GouveiaTorresCMAM} and Example~\ref{ex:prob2358}), with
different dependences of the infinitesimal generators:
\small{\emph{D1} -- $[T(t,\mathbf{x},\mathbf{u},\boldmath{\psi})$,
$\mathbf{X}(t,\mathbf{x},\mathbf{u},\boldmath{\psi})$,
$\mathbf{U}(t,\mathbf{x},\mathbf{u},\boldmath{\psi})$,
$\boldsymbol{\Psi}(t,\mathbf{x},\mathbf{u},\boldmath{\psi})]$;
\emph{D2} -- $[T(t)$, $\mathbf{X}(t,\mathbf{x})$,
$\mathbf{U}(\mathbf{u},\boldmath{\psi})$,
$\boldsymbol{\Psi}(\boldmath{\psi})]$; \emph{D3} -- $[T(t)$,
$\mathbf{X}(t,\mathbf{x})$, $\mathbf{U}(t,\mathbf{u})$,
$\boldsymbol{\Psi}(t,\boldmath{\psi})]$; \emph{D4} -- $[T(t)$,
$\mathbf{X}(\mathbf{x})$, $\mathbf{U}(\mathbf{u})$,
$\boldsymbol{\Psi}(\boldmath{\psi})]$}. } \label{tab:depGerad}
\end{center}
\end{table}

We verify that of the four sets of studied generators, just with
\emph{D4} it was not possible to obtain, with full generality, the
totality of Noether's conservation laws for the twenty two
considered problems. The set of generators \emph{D3} ($T(t) $,
$\mathbf{X}(t,\mathbf{x}) $, $\mathbf{U}(t,\mathbf{u}) $,
$\boldsymbol{\Psi}(t,\boldsymbol{\psi}) $) gives the best
compromise: it presents the best running times, between the
generators that give the complete family of variational symmetries
and Noether conservation laws for the problems we have studied;
running times are much better than the ones obtained with the
generators \emph{D1}. We recommend the user to try configuration
\emph{D3} first on his/her own optimal control problems.
Considering $t$ and $\mathbf{x}$ for the dependences of the gauge
term -- $G(t,\mathbf{x})$ -- the system of PDEs that we have to
solve, in order to find the variational symmetries, takes form
(\textrm{cf.} \eqref{eq:detGerad3})
\begin{equation}
\label{eq:detGerad3a}
\left\{
\begin{array}{l}
\displaystyle \medskip
\frac{\partial H}{\partial t}T
+\frac{\partial H}{\partial \mathbf{x}}\cdot \mathbf{X}
+\frac{\partial H}{\partial \mathbf{u}}\cdot \mathbf{U}
+\frac{\partial H}{\partial \boldsymbol{\psi}}\cdot \boldsymbol{\Psi}
+H \frac{\partial T}{\partial t}
-\boldsymbol{\psi}\T \cdot
\frac{\partial \mathbf{X}}{\partial t}
-\frac{\partial G}{\partial t}=0 \, ,\\
\displaystyle \medskip \boldsymbol{\Psi}\T +\boldsymbol{\psi}\T
\cdot \frac{\partial \mathbf{X}} {\partial \mathbf{x}} +
\frac{\partial G}{\partial \mathbf{x}} =\mathbf{0} \, .
\end{array}
\right.
\end{equation}
Our present procedure \emph{Symmetry} computes, by default, the
variational symmetries as defined by \emph{D3}, and with a gauge
term $G(t,\mathbf{x})$: by default \emph{Symmetry} solves system
\eqref{eq:detGerad3a}. Through optional parameters, it is possible
to find the variational symmetries for other generators and gauge
terms: in order to use all the dependences (\emph{D1}) one must
use option~\texttt{alldep}; to use a minimum of dependences
(\emph{D4}) one uses option~\texttt{mindep}. We remark that with
the class of generators \emph{D3}, $T$ and $\mathbf{X}$ are not
functions of $\mathbf{u}$, and there is no need to assume the
control variables $\mathbf{u}$ to be smooth functions
(\textrm{cf.} \eqref{eq:smothCtrl}).

Table~\ref{tab:depGerad2} shows the computing running times needed
to obtain all the variational symmetries of the problems in
\cite[\S4, \S5]{GouveiaTorresCMAM}, by using the default version
of procedure \emph{Symmetry} we give here (generators \emph{D3});
and by using the version in \cite{GouveiaTorresCMAM}, which is a
particular case of our present procedure -- see \S\ref{sec:exOC}
for examples not covered by the previous methods in
\cite{GouveiaTorresCMAM} -- obtained using option \texttt{alldep},
that is, generators \emph{D1}. The time needed to compute the
variational symmetries for the $(2,3,5)$ problem (Example~4.6 in
\cite{GouveiaTorresCMAM}) decreased from thirty minutes to one.

\begin{table}[!h]
\begin{center}
\small
\begin{tabular}{|c|r|r|r|r|r|r|r|r|r|r|}
\hline
 &4.1&4.2&4.3& 4.4& 4.5 & 4.6
& 5.1& 5.2& 5.3& 5.4 \\
\hline \hline
\emph{D1} & $2''$  & $1'13''$ & $2'44''$ &  $6'41''$ & $1'04''$
& $30'34''$  & $8''$ & $17''$ &  $6'42''$ & $1''$\\
\hline
\emph{D3} & $0''$  & $5''$ & $11''$ &  $18''$ & $4''$
& $1'09''$  & $0''$ & $3''$ &  $16''$ & $0''$\\
\hline
\end{tabular}
\normalsize
\caption{Running times of procedure \emph{Symmetry} for all the
problems of previous work \cite{GouveiaTorresCMAM}, with the
generator sets \emph{D1} (the only possibility in
\cite{GouveiaTorresCMAM}) --
\small{$[T(t,\mathbf{x},\mathbf{u},\boldmath{\psi})$,
$\mathbf{X}(t,\mathbf{x},\mathbf{u},\boldmath{\psi})$,
$\mathbf{U}(t,\mathbf{x},\mathbf{u},\boldmath{\psi})$,
$\boldsymbol{\Psi}(t,\mathbf{x},\mathbf{u},\boldmath{\psi})]$},
and \emph{D3} -- \small{$[T(t)$, $\mathbf{X}(t,\mathbf{x})$,
$\mathbf{U}(t,\mathbf{u})$,
$\boldsymbol{\Psi}(t,\boldmath{\psi})]$}.} \label{tab:depGerad2}
\end{center}
\end{table}

The use of generators with a smaller number of dependences leads
to a drastic reduction of the computing running times. For the
studied problems, the use of generators \emph{D3} permit to obtain
the same results while decrease the total processing times for
about 4\% of the ones verified in \cite{GouveiaTorresCMAM}
(generators \emph{D1}).


\section{Examples of the new possibilities}
\label{sec:exOC}

In order to show the functionality and the use of the new
procedures, we apply our \maple package to three concrete
optimal control problems which are not covered by the previous
results in \cite{gouv04,GouveiaTorresCMAM}. All the examples were
carried out with \maple version 10 on a 1.4GHz 512MB RAM Pentium
Centrino. The running time of procedure \emph{Symmetry} is
indicated, for each example, in the format min'sec''. All the
other \maple commands run instantaneously.


\subsection{Variational symmetries up to a gauge term}
\label{sec:exCV}

We begin with a very simple example
of the classical calculus of variations.
We recall that for the fundamental problem of the calculus of variations
there are no abnormal extremals, so one can choose $\psi_0 = -1$
(we use option \texttt{noabn} of our \maple package).


\begin{example}
\label{ex:Ex1_1} (0'00'')
Let us consider the following scalar problem of the calculus of variations
($n = m = 1$):
\begin{gather*}
\int_a^b \left(u(t)\right)^2\,\dd t \longrightarrow \min \, , \\
\dot{x}(t) = u(t)\, .
\end{gather*}
In this case $L = u^2$ and $\varphi = u$.
First we obtain the variational symmetries of the problem
(\maple procedure \emph{Symmetry})
up to a gauge term (parameter \texttt{gauge}).
\small
\begin{verbatim}
> S := Symmetry(u^2,u,t,x,u,showt,gauge);
\end{verbatim}
\begin{eqnarray*}
 S:=
\biggl[T=2\,C_{{2}}t+C_{{6}},\; X=\frac{1}{2}\,{\frac {C_{{3}}t}{\psi_{{0}}}}
+C_{{2}}x
(t) +C_{{4}},\; U=\frac{1}{2}\,{\frac{C_{{3}}}{\psi_{{0}}}}-u(t) C_{{2}}
,\\
\Psi=-\psi(t) C_{{2}}-C_{{3}},\;
 {\it GAUGE}=C_{{3}}x
(t) +C_{{5}}\biggr]
\end{eqnarray*}
\normalsize Noether conservation laws are obtained through
Theorem~\ref{thm:TNoether2} (\maple procedure \emph{Noether}) with
the generators and the gauge term just obtained. \small
\begin{verbatim}
> CL := Noether(u^2,u,t,x,u,S,showt,noabn,H);
\end{verbatim}
\[
 CL:=\left( -\frac{1}{2}\,C_{{3}}t+C_{{2}}x(t) +C_{{4}} \right) \psi
(t) -H \left( 2\,C_{{2}}t+C_{{6}} \right) +C_{{3}}x
(t) +C_{{5}}={\it const}
\]
\normalsize
The Hamiltonian $H$, which appears in the above family of conservation laws,
is given by \eqref{eq:hamilt}:
\small
\begin{verbatim}
> H := PMP(u^2,u,t,x,u, evalH,showt,noabn);
\end{verbatim}
\[
H :=- u(t)^{2}+u(t) \psi
(t)
\]
\normalsize
This is a very simple problem, just used to illustrate,
in the simplest possible way, our \maple procedures. In this case
it is an easy exercise to obtain the extremals by direct application
of the Pontryagin Maximum Principle or the Euler-Lagrange equations,
\small
\begin{verbatim}
> extremals := PMP(u^2,u,t,x,u,showt,noabn);
\end{verbatim}
\[
{\it extremals}:= \left\{ \psi(t) =K_2,\;
x(t) =\frac{1}{2}\,K_2\,t+K_1,\; u(t)
 =\frac{1}{2}\,K_2 \right\}
\]
\normalsize
and one can validate the obtained conservation laws
by applying the definition of conservation law: by definition,
the obtained family of conservation laws must hold
along all the extremals of the problem.
\small
\begin{verbatim}
> subs(extremals,CL);
\end{verbatim}
\[
K_2\,C_{{2}}K_1+K_2\,C_{{4}}-\frac{1}{4}\,K_2^{2}
C_{{6}}+C_{{3}}K_1+C_{{5}}={\it const}
\]
\normalsize
\end{example}


\subsection{Presence of nonconservative forces}

We consider now a problem of the calculus of variations
under the action of a nonconservative force.
The problem is borrowed from \cite[\S 4]{CD:Djukic:1980}.


\begin{example}
\label{ex:fnc3} ($n=1$, $m=2$, 0'01'') The problem is defined by the
Lagrangian $L(q,\dot q, \ddot q)=\frac{1}{2}\,{\ddot
q(t)}^{2}+\frac{1}{2}\,a{\dot q(t)}^{2}+\frac{1}{2}\,b{q(t)}^{2}$,
and presence of the nonconservative force $f(t)=\mu \,\dot q(t)+{\frac {{\mu}^{2}
}{{a}^{2}}}\ddot q(t)-2\ {\frac {\mu}{a}} \dddot{q}(t)$ which depends
on higher-order derivatives ($a$, $b$, and $\mu$ are constants).
\small
\begin{verbatim}
> PDEtools[declare](prime=t);
\end{verbatim}
\textit{derivatives with respect to t of functions of one variable
will now be displayed with '}
\begin{verbatim}
> L := u^2/2+a*v^2/2+b*q^2/2;
> phi := [v,u];
> f := mu*v+mu^2/a^2*u-2*mu/a*z(t);
\end{verbatim}
\begin{eqnarray*}
L&:=&\frac{1}{2}\,{u}^{2}+\frac{1}{2}\,a{v}^{2}+\frac{1}{2}\,b{q}^{2}\\
\varphi &:=&[v,\; u]\\
f&:=&\mu\,v+{\frac {{\mu}^{2}u}{{a}^{2}}}-2\,{\frac {\mu\,z(t) }{a}}
\end{eqnarray*}
\begin{verbatim}
> S := Symmetry(L, phi, t, [q,v], u);
\end{verbatim}
\[
S := [T=C_{{1}},\; X_{{1}}=0,\; X_{{2}}=0,\; U=0,\; \Psi_{{1}}=0,\; \Psi_{{2}}=0]
\]
\normalsize
\small
\begin{verbatim}
> CL := Noether(L, phi, t, [q,v], u, S, ncf=[f,0], noabn);
\end{verbatim}
\begin{eqnarray*}
CL:=- \left( -\frac{1}{2}\, u(t)^{2}-\frac{1}{2}\,a \,
v(t) ^{2}-\frac{1}{2}\,b \,q(t)^{2}+\psi_{{1}}(t) v(t) +\psi_{{2
}}(t) u(t)  \right) C_{{1}}\\
+\int \!C_{{1}}
{\it q'}\, \left( \mu\,v(t) +{\frac {{\mu}^{2}u(t)
 }{{a}^{2}}}-2\,{\frac {\mu\,z(t) }{a}} \right)
{dt}={\it const}\end{eqnarray*}
\normalsize
The multipliers $\psi_1(t)$ and $\psi_2(t)$ are obtained using
the adjoint system and the stationary condition, as given
by Theorem~\ref{thm:pMaxPont}.
\small
\begin{verbatim}
> sys := PMP(L, phi, t, [q,v], u, noabn, evalSyst, ncf=[f,0], showt);
\end{verbatim}
\begin{eqnarray*}
{\it sys}:=\biggl[ \left\{ {\it q'}=v(t) ,{\it v'}=u(t)
  \right\} , \biggl\{ -{\psi_{{1}}}^{'}=-\mu\,v(t)
   -{\frac {{\mu}^{2}u(t) }{{a}^{2}}}+2\,{\frac {
\mu\,z(t) }{a}}-bq(t) , \\
-{\psi_{{2}}}^{'}=-av(t) +
\psi_{{1}}(t) \biggr\} , \left\{
-u(t) +\psi_{{2}}(t) =0 \right\} \biggr]
\end{eqnarray*}
\normalsize
\small
\begin{verbatim}
> dsolve({sys[2][2],sys[3][]},{psi[1](t),psi[2](t)});
\end{verbatim}
\[
\left\{ \psi_{{2}}(t) =u(t) ,\; \psi_{{1}}(t) =-{\it u'}+av(t)  \right\}
\]
\normalsize
With substitutions
\small
\begin{verbatim}
> subs(%, z(t)=diff(u(t),t), u(t)=diff(v(t),t), v(t)=diff(q(t),t),
                                                             C[1]=1, CL);
\end{verbatim}
\[
-\frac{1}{2}\,{{\it q''}}^{2}+\frac{1}{2}\,a{{\it q'}}^{2}+\frac{1}{2}\,b \,
 q(t)^{2}
- \left( -{\it q'''}+a{\it q'} \right) {\it q'}+
\int \!{\it q'}\, \left( \mu\,{\it q'}+{\frac {{\mu}^{2}{\it q''}}{{a}
^{2}}}-2\,{\frac {\mu\,{\it q'''}}{a}} \right) {dt}={\it const}
\]
\normalsize
one obtains the conservation law \cite[\S 4]{CD:Djukic:1980}.
We remark that the conclusion is nontrivial, and difficult to obtain
without Noether's principle.
\end{example}


\subsection{The sub-Riemannian nilpotent case $(2,3,5,8)$}

We finish the section by applying our \maple package to one
important problem: the study of sub-Riemannian geodesics. The
reader, interested in the study of symmetries of flat
distributions of sub-Riemannian geometry, is referred to
\cite{Sachkov04}. Here we use a formulation of the nilpotent
problem $(2,3,5,8)$ which is obtained using the results of
\cite{RochaPhD}.


\begin{example}
\label{ex:prob2358} (44'16'') The problem can be defined in the
following way:
\begin{equation*}
\frac{1}{2}\,\int_a^b \left({u_1(t)}^2+{u_2(t)}^{2}\right)\dd t
\longrightarrow \min \, , \quad
\begin{cases}
\dot x_1(t)=u_1(t)\, ,\\
\dot x_2(t)=u_2(t)\, ,\\
\dot x_3(t)=u_2(t) x_1(t)\, ,\\
\dot x_4(t)=\frac{1}{2}\,u_2(t) {x_1(t)}^2\, ,\\
\dot x_5(t)=u_2(t) x_1(t) x_2(t)\, ,\\
\dot x_6(t)=\frac{1}{6}\,u_2(t) {x_1(t)}^3\, ,\\
\dot x_7(t)=\frac{1}{2}\,u_2(t) {x_1(t)}^2 x_2(t)\, ,\\
\dot x_8(t)=\frac{1}{2}\,u_2(t) {x_1(t)} {x_2(t)}^2\, .\\
\end{cases}
\end{equation*}
The integrability of the problem is still an open question
\cite{RochaTorres4thJM,Sachkov04}, but eight independent conservation laws
can be determined with our present \maple package.
\small
\begin{verbatim}
> L := 1/2*(u[1]^2+u[2]^2);
> phi:=[u[1], u[2], u[2]*x[1], (u[2]/2)*x[1]^2, u[2]*x[1]*x[2],
            (u[2]/6)*x[1]^3, (u[2]/2)*x[1]^2*x[2], (u[2]/2)*x[1]*x[2]^2];
> XX := [x[i]$i=1..8];
> UU := [u[1],u[2]];
\end{verbatim}
\vspace*{-0.3cm}
\begin{eqnarray*}
L&:=&\frac{1}{2}\,{u_1}^2+\frac{1}{2}\,{u_2}^2\\
\varphi &:=&\left[u_{{1}},\; u_{{2}},\; u_{{2}}x_{{1}},\frac{1}{2}\,
u_{{2}}{x_{{1}}}^{2},\; u_{{2}}x_{{1
}}x_{{2}},\frac{1}{6}\, u_{{2}}{x_{{1}}}^{3},\frac{1}{2}\, u_{{2}}
{x_{{1}}}^{2}x_{{2}},\frac{1}{2}\,
u_{{2}}x_{{1}}{x_{{2}}}^{2}\right]\\
XX&:=&[x_{{1}},x_{{2}},x_{{3}},x_{{4}},x_{{5}},x_{{6}},x_{{7}},x_{{8}}]\\
UU&:=&[u_{{1}},\; u_{{2}}]
\end{eqnarray*}
\normalsize
\small
\begin{verbatim}
> Symmetry(L, phi, t, XX, UU);
\end{verbatim}
\vspace*{-0.2cm}
\begin{multline*}
\biggl[T=C_{{1}}t+C_{{7}},\; X_{{1}}=\frac{1}{2}\,C_{{1}}x_{{1}},\;
X_{{2}}=C_{{2}}+\frac{1}{2}\,C
_{{1}}x_{{2}},\; X_{{3}}=C_{{1}}x_{{3}}+C_{{8}},\;\\
X_{{4}}=\frac{3}{2}\,
C_{{1}}x_{{4
}}+C_{{6}},\;
X_{{5}}=C_{{2}}x_{{3}}+\frac{3}{2}\,C_{{1}}x_{{5}}+C_{{3}},\; X_{{6}}=
2\,C_{{1}}x_{{6}}+C_{{5}},\; \\
X_{{7}}=C_{{2}}x_{{4}}+2\,C_{{1}}x_{{7}}+C_{
{9}},\;
X_{{8}}=C_{{2}}x_{{5}}+2\,C_{{1}}x_{{8}}+C_{{4}},\; U_{{1}}=-\frac{1}{2}\,u_
{{1}}C_{{1}},\;\\
U_{{2}}=-\frac{1}{2}\,C_{{1}}u_{{2}},\; \Psi_{{1}}=-\frac{1}{2}
\,C_{{1}}\psi
_{{1}},
\Psi_{{2}}=-\frac{1}{2}\,C_{{1}}\psi_{{2}},\; \Psi_{{3}}=-\psi_{{3}}C_{{1}
}-C_{{2}}\psi_{{5}},\;\\
\Psi_{{4}}=-\frac{3}{2}\,\psi_{{4}}C_{{1}}-C_{{2}}
\psi_{{7
}},\;
\Psi_{{5}}=-\frac{3}{2}\,C_{{1}}\psi_{{5}}-C_{{2}}\psi_{{8}},\; \Psi_{{6}}=-2
\,C_{{1}}\psi_{{6}},\;\\
\Psi_{{7}}=-2\,C_{{1}}\psi_{{7}},\; \Psi_{{8}}=-2\,C_
{{1}}\psi_{{8}}\biggr]
\end{multline*}
\normalsize
\small
\begin{verbatim}
> CL := Noether(L, phi, t, XX, UU, %, H);
\end{verbatim}
\begin{multline*}
 CL:=\frac{1}{2}\,C_{{1}}x_{{1}}\psi_{{1}}+ \left( C_{{2}}+\frac{1}{2}\,
 C_{{1}}x_{{2}}
 \right) \psi_{{2}}+ \left( C_{{1}}x_{{3}}+C_{{8}} \right) \psi_{{3}}+
 \left( \frac{3}{2}\,C_{{1}}x_{{4}}+C_{{6}} \right) \psi_{{4}}\\
 + \left( C_{{2}}
x_{{3}}+\frac{3}{2}\,C_{{1}}x_{{5}}+C_{{3}} \right) \psi_{{5}}+ \left( 2\,C_{{
1}}x_{{6}}+C_{{5}} \right) \psi_{{6}}+ \left( C_{{2}}x_{{4}}+2\,C_{{1}
}x_{{7}}+C_{{9}} \right) \psi_{{7}}\\
+ \left( C_{{2}}x_{{5}}+2\,C_{{1}}x
_{{8}}+C_{{4}} \right) \psi_{{8}}-H \left( C_{{1}}t+C_{{7}} \right) ={
\it const}
\end{multline*}
\normalsize
The Hamiltonian is given by
\small
\begin{verbatim}
> Hamilt := PMP(L, phi, t, XX, UU, noabn, evalH);
\end{verbatim}
\begin{multline*}
{\it Hamilt}:=-\frac{1}{2}\,{u_{{1}}}^{2}-\frac{1}{2}\,{u_{{2}}}^{2}
+\psi_{{1}}u_{{1}}+\psi_{{2}}u_{
{2}}+\psi_{{3}}u_{{2}}x_{{1}}+\frac{1}{2}\,\psi_{{4}}u_{{2}}{x_{{1}}}^{2}+\psi
_{{5}}u_{{2}}x_{{1}}x_{{2}}\\+\frac{1}{6}\,u_{{2}}{x_{{1}}}^{3}\psi_{{6}}
+\frac{1}{2}\,u
_{{2}}{x_{{1}}}^{2}x_{{2}}\psi_{{7}}+\frac{1}{2}\,u_{{2}}x_{{1}}{x_{{2}}}^{2}
\psi_{{8}}
\end{multline*}
\normalsize
and the extremal controls are obtained through the stationary condition.
\small
\begin{verbatim}
> PMP(L,phi,t, XX, UU, noabn, evalSyst)[3];
\end{verbatim}
\begin{eqnarray*}
\biggl\{
-u_{{2}}+\psi_{{2}}+\psi_{{3}}x_{{1}}
+\frac{1}{2}\,\psi_{{4}}{x_{{1}}}^{2}
+\psi_{{5}}x_{{1}}x_{{2}}+\frac{1}{6}\,
{x_{{1}}}^{3}\psi_{{6}}+\frac{1}{2}\,{x_{{1}}}^{2}x_{{2}}\psi_{{7}}
+\frac{1}{2}\,
x_{{1}}{x_{{2}}}^{2}\psi_{{8}}=0,\\
-u_{{1}}+\psi_{{1}}=0
\biggr\}
\end{eqnarray*}
\normalsize
\small
\begin{verbatim}
> solve(%,{u[1],u[2]});
\end{verbatim}
\[
\left\{ u_{{1}}=\psi_{{1}}, \,u_{{2}}=\psi_{{5}}x_{{1}}x_{{2}}+\psi_{{2}
}+\psi_{{3}}x_{{1}}+\frac{1}{2}\,\psi_{{4}}{x_{{1}}}^{2}+\frac{1}{6}\,
{x_{{1}}}^{3}
\psi_{{6}}+\frac{1}{2}\,{x_{{1}}}^{2}x_{{2}}\psi_{{7}}+\frac{1}{2}\,
x_{{1}}{x_{{2}}}^{
2}\psi_{{8}} \right\}
\]
\normalsize
\small
\begin{verbatim}
> H = expand(subs(%, Hamilt));
\end{verbatim}
\begin{eqnarray*}
H=\frac{1}{2}\,\psi_{{2}}x_{{1}}{x_{{2}}}^{2}\psi_{{8}}+\psi_{{5}}x_{{1}}x_{{2
}}\psi_{{2}}+\psi_{{5}}{x_{{1}}}^{2}x_{{2}}\psi_{{3}}+\frac{1}{2}\,\psi_{{2}}
\psi_{{4}}{x_{{1}}}^{2}+\frac{1}{2}\,\psi_{{3}}{x_{{1}}}^{3}\psi_{{4}}
\qquad\qquad\qquad\\
+\frac{1}{2}\,{\psi_{{5}}}^{2}{x_{{1}}}^{2}{x_{{2}}}^{2}
+\frac{1}{6}\,\psi_{{2}}{x_{{1}}}^{3}
\psi_{{6}}+\frac{1}{8}\,{x_{{1}}}^{2}{x_{{2}}}^{4}{\psi_{{8}}}^{2}+
\frac{1}{8}\,{x_{{1
}}}^{4}{x_{{2}}}^{2}{\psi_{{7}}}^{2}+\frac{1}{12}\,\psi_{{4}}{x_{{1}}}^{5}\psi
_{{6}}\\
+\frac{1}{2}\,{\psi_{{3}}}^{2}{x_{{1}}}^{2}
+\frac{1}{72}\,{x_{{1}}}^{
6}{\psi_{{6}}}^{2}+\frac{1}{2}\,{\psi_{{2}}}^{2}
+\frac{1}{2}\,{\psi_{{1}}}^{2}+\frac{1}{8}\,{
\psi_{{4}}}^{2}{x_{{1}}}^{4}+\frac{1}{6}\,\psi_{{5}}{x_{{1}}}^{4}x_{{2}}\psi_{
{6}}+\frac{1}{2}\,\psi_{{3}}{x_{{1}}}^{2}{x_{{2}}}^{2}\psi_{{8}}\\
+\frac{1}{4}\,\psi_{{4}}{x_{{1}}}^{3}{x_{{2}}}^{2}\psi_{{8}}
+\frac{1}{4}\,\psi_{{4}}{x_{{1}}}^{4}x_{
{2}}\psi_{{7}}+\psi_{{2}}\psi_{{3}}x_{{1}}
+\frac{1}{4}\,{x_{{1}}}^{3}{x_{{2}}}
^{3}\psi_{{7}}\psi_{{8}}+\frac{1}{12}\,{x_{{1}}}^{5}\psi_{{6}}x_{{2}}\psi_{{7}
}\\
+\frac{1}{12}\,{x_{{1}}}^{4}\psi_{{6}}{x_{{2}}}^{2}\psi_{{8}}
+\frac{1}{2}\,\psi_{{2}}
{x_{{1}}}^{2}x_{{2}}\psi_{{7}}+\frac{1}{2}\,\psi_{{5}}{x_{{1}}}^{3}x_{{2}}\psi
_{{4}}
+\frac{1}{2}\,\psi_{{5}}{x_{{1}}}^{2}{x_{{2}}}^{3}\psi_{{8}}+\frac{1}{2}\,
\psi_{
{5}}{x_{{1}}}^{3}{x_{{2}}}^{2}\psi_{{7}}\\
+\frac{1}{2}\,\psi_{{3}}{x_{{1}}}^{3}x
_{{2}}\psi_{{7}}+\frac{1}{6}\,\psi_{{3}}{x_{{1}}}^{4}\psi_{{6}}
\end{eqnarray*}
\normalsize
Now, the eight conservation laws, we are looking for, are easily obtained:
\small
\begin{verbatim}
> subs(C[8]= 1, seq(C[i]=0,i=1..9), CL);
> subs(C[6]= 1, seq(C[i]=0,i=1..9), CL);
> subs(C[3]= 1, seq(C[i]=0,i=1..9), CL);
> subs(C[5]= 1, seq(C[i]=0,i=1..9), CL);
> subs(C[9]= 1, seq(C[i]=0,i=1..9), CL);
> subs(C[4]= 1, seq(C[i]=0,i=1..9), CL);
> subs(C[2]= 1, seq(C[i]=0,i=1..9), CL);
> subs(C[7]=-1, seq(C[i]=0,i=1..9), CL);
\end{verbatim}
\begin{eqnarray*}
&&\psi_{{3}} ={\it const}\\
&&\psi_{{4}} ={\it const}\\
&&\psi_{{5}} ={\it const}\\
&&\psi_{{6}} ={\it const}\\
&&\psi_{{7}} ={\it const}\\
&&\psi_{{8}} ={\it const}\\
&&\psi_{{2}}+x_{{3}}\psi_{{5}}+x_{{4}}\psi_{{7}}+x_{{5}}\psi_{{8}}={\it
const}\\
&&H={\it const}
\end{eqnarray*}
\normalsize
\end{example}

Given the results of \cite{RochaPhD}, one can say that the
sub-Riemannian nilpotent Lie group of type $(2,3,5,8)$ has seven
trivial first integrals: the Hamiltonian $H$; and the multipliers
$\psi_3$, $\psi_4$, $\psi_5$, $\psi_6$, $\psi_7$, $\psi_8$.
Together with the non-trivial first integral
$\psi_{{2}}+x_{{3}}\psi_{{5}}+x_{{4}}\psi_{{7}}+x_{{5}}\psi_{{8}}$,
here first obtained, it is possible to prove that the system is
integrable. This is nontrivial since Liouville theorem does not
apply: the set of first integrals is not involutive (for instance,
Poisson bracket between $\psi_{{3}}$ and
$\psi_{{2}}+x_{{3}}\psi_{{5}}+x_{{4}}\psi_{{7}}+x_{{5}}\psi_{{8}}$
is not zero). This question is under study and will be addressed
in a forthcoming publication.


\section*{Acknowledgements}

The support from the control theory group
(\textsf{cotg}) of the R\&D unit \textsf{CEOC}
is here acknowledged. PG was also supported by the program
PRODEP III/5.3/2003; DT and ER by the research project
``Advances in Nonlinear Control and Calculus of
Variations'' POCTI/MAT/41683/2001.



\begin{thebibliography}{99}

\bibitem{cotcot} B. Bonnard, J.-B. Caillau, E. Trélat.
\textit{Cotcot: short reference manual},
Ecole Nationale Sup\'{e}rieure
d'Electronique, d'Electrotechnique
d'Informatique, d'Hydraulique et de T\'{e}l\'{e}com,
Institut de Recherche en Informatique de Toulouse,
Technical Report RT/APO/05/1.

\bibitem{Terrab95} E. S. Cheb-Terrab, K. von Bulow.
\textit{A Computational approach for the analytical solving of
partial differential equations}, Computer Physics Communications,
{\bf 90} (1995), pp.~102--116.

\bibitem{Djukic73} D. S. Djukic.
\textit{Noether's theorem for optimum control systems}, Internat.
J. Control, {\bf 1} (1973), no.~18, pp.~667--672.

\bibitem{CD:Djukic:1980} D. S. Djukic, A. M. Strauss.
\textit{Noether's theory for nonconservative generalised mechanical systems},
J. Phys. A {\bf 13} (1980), no.~2, 431--435.

\bibitem{gast05} Gast\~{a}o S. F. Frederico, Delfim F. M. Torres.
\textit{Nonconservative Noether's Theorem in Optimal Control},
Proc. 13th IFAC Workshop on ``Control Applications of
Optimisation'' (CAO'06), 26-28 April 2006, Paris - Cachan, France.
IFAC publication, Elsevier Ltd, Oxford, UK (in press).

\bibitem{Fu03} Jing-Li Fu, Li-Qun Chen.
\textit{Non-Noether symmetries and conserved quantities of
nonconservative dynamical systems}, Phys. Lett. A {\bf 317}
(2003), no.~3-4, 255--259.

\bibitem{gouv04}  Paulo D. F. Gouveia, Delfim F. M. Torres.
\textit{Computaç\~{a}o Alg\'{e}brica no C\'{a}lculo das
Variaç\~{o}es: Determinaç\~{a}o de Simetrias e Leis de
Conservaç\~{a}o} (in Portuguese), TEMA Tend. Mat. Apl. Comput.
{\bf 6} (2005), no.~1, pp.~81--90.

\bibitem{GouveiaTorresCMAM} Paulo D. F. Gouveia, Delfim F. M. Torres.
\textit{Automatic Computation of Conservation Laws in the Calculus
of Variations and Optimal Control}, Comput. Methods Appl. Math.
{\bf 5} (2005), no.~4, pp.~387--409.

\bibitem{Noether} E. Noether.
\textit{Invariante Variationsprobleme},
G\"{o}tt. Nachr. (1918), pp.~235--257.

\bibitem{Pontryagin62} L. S. Pontryagin, V. G. Boltyanskii,
R. V. Gamkrelidze, E. F. Mishchenko. \textit{The mathematical
theory of optimal processes}, Interscience Publishers John Wiley
\& Sons, Inc. New York-London, 1962.

\bibitem{RochaPhD} Eug\'{e}nio A. M. Rocha.
\textit{An Algebraic Approach to Nonlinear Control Theory}, PhD
thesis, University of Aveiro, 2004.

\bibitem{RochaTorres4thJM} Eug\'enio A. M. Rocha, Delfim F. M. Torres.
\textit{Quadratures of Pontryagin Extremals for Optimal Control
Problems}, Proceedings of the 4th Junior European Meeting on
``Control and Optimization'', Bia\l ystok Technical University,
Bia\l ystok, Poland, 11--14 September 2005, Control \& Cybernetics
(accepted).

\bibitem{Sachkov04} Yu. L. Sachkov.
\textit{Symmetries of flat rank two distributions and sub-Riemannian structures},
Trans. Amer. Math. Soc., {\bf 356} (2004), no.~2, pp.~457--494.

\bibitem{delfimEJC} Delfim F. M. Torres.
\textit{On the Noether Theorem for Optimal Control},
European Journal of Control, {\bf 8} (2002), no.~1, pp.~56--63.

\bibitem{Torres04} Delfim F. M. Torres.
\textit{Quasi-Invariant Optimal Control Problems},
Portugali\ae\ Mathematica (N.S.), {\bf 61} (2004),
no.~1, pp.~97--114.

\bibitem{Torres05} Delfim F. M. Torres.
\textit{Weak Conservation Laws for Minimizers which are not
Pontryagin Extremals}, Proc. of the 2005 International Conference
``Physics and Control'' (PhysCon 2005), August 24-26, 2005, Saint
Petersburg, Russia.  Edited by A.L.~Fradkov and A.N.~Churilov,
2005 IEEE, pp.~134--138.

\end{thebibliography}
\end{document}